\documentclass{amsart}

\newtheorem{theorem}{Theorem}[section]
\newtheorem{lemma}[theorem]{Lemma}
\newtheorem{corollary}[theorem]{Corollary}
\newtheorem{proposition}[theorem]{Proposition}
\newtheorem{remark}[theorem]{Remark}

\theoremstyle{definition}

\theoremstyle{remark}

\numberwithin{equation}{section}

\begin{document}
\title[Power Moments of Ternary Kloosterman Sums]{$\begin{array}{c}
         \text{Infinite Families of Recursive Formulas Generating}\\
         \text{Power Moments of Ternary Kloosterman Sums with}\\
           \text{ Trace Nonzero Square Arguments:
           $O(2n+1,2^{r})$ Case}
       \end{array}$}

\author{dae san kim}
\address{Department of Mathematics, Sogang University, Seoul 121-742, Korea}
\curraddr{Department of Mathematics, Sogang University, Seoul
121-742, Korea} \email{dskim@sogong.ac.kr}
\thanks{}

\subjclass[]{}

\date{}

\dedicatory{ }

\keywords{}

\begin{abstract}
In this paper, we construct four infinite families  of ternary
linear codes  associated with double cosets in $O(2n+1,q)$ with
respect to certain maximal parabolic subgroup  of the special
orthogonal  group $SO(2n+1,q)$. Here $q$ is a power of three. Then
we obtain two infinite families of recursive formulas, the one
generating the power moments of Kloosterman sums with $\lq\lq$trace
nonzero square arguments" and the other generating the even power
moments of those. Both of these families are expressed in terms of
the frequencies of weights in the codes associated with those double
cosets in $O(2n+1,q)$ and in the codes associated with similar
double cosets in the symplectic group $Sp(2n,q)$. This is done via
Pless power moment identity and by utilizing the explicit
expressions of exponential sums over those double cosets related to
the evaluations of $\lq\lq$Gauss sums" for the orthogonal group
$O(2n+1,q)$.\\

  Index terms- power moment,  Kloosterman sum,  trace  nonzero square
argument, orthogonal group, symplectic group, double cosets, maximal
parabolic subgroup, Pless power moment identity, weight
distribution,
Gauss sum.\\

MSC 2000: 11T23, 20G40, 94B05.

\end{abstract}

\maketitle

\section{Introduction}
 Let $\psi$ be a nontrivial additive character of the finite field $\mathbb {F}_{q}$ with $q=p^{r}$ elements
 ($p$ a prime). Then the Kloosterman sum $K(\psi;a)$(\cite{LN1}) is defined by
 \begin{align*}
K(\psi;a)=\sum_{\alpha\in\mathbb{F}_{q}^{*}}\psi(\alpha_{}^{}+a\alpha_{}^{-1})
~(a\in\mathbb{F}_{q}^{*}).
\end{align*}

The Kloosterman sum was introduced in 1926(\cite{K1}) to give an
estimate for the Fourier coefficients of modular forms.

For each nonnegative integer $h$, by $MK(\psi)^{h}$ we will denote
the $h$-th moment of the Kloosterman sum $K(\psi;a)$. Namely, it is
given by
\begin{align*}
MK(\psi)^{h}=\sum_{a\in\mathbb{F}_{q}^{*}}K(\psi;a)^{h}.
\end{align*}
If $\psi=\lambda$ is the canonical additive character of
$\mathbb{F}_{q}$, then $MK(\lambda)^{h}$ will be simply denoted by
$MK^{h}$.

  Explicit computations on power moments of Kloosterman sums were
begun with the paper \cite{S1} of Sali\'{e} in 1931, where he
showed, for any odd prime $q$,
\begin{align*}
MK^h=q^2 M_{h-1} -(q-1)^{h-1}+2(-1)^{h-1}~(h\geq1).
\end{align*}
Here $M_0=0$, and, for $h\in\mathbb{Z}_{>0}$,
\begin{align*}
M_h=|\{(\alpha_1,\cdots,\alpha_h)\in(\mathbb{F}_q^*)^h|\sum_{j=1}^h
\alpha_j=1=\sum_{j=1}^h \alpha_j^{-1}\}|.
\end{align*}
For $q=p$ odd prime, Sali\'{e} obtained $MK^1$, $MK^2$, $MK^3$,
$MK^4$ in \cite{S1} by determining $M_1$, $M_2$, $M_3$. On the other
hand, $MK^5$ can be expressed in terms of the $p$-th eigenvalue for
a weight 3 newform on $\Gamma_0(15)$(cf. \cite{L1}, \cite{PV1}).
$MK^6$ can be expressed in terms of the $p$-th eigenvalue for a
weight 4 newform on $\Gamma_0(6)$ (cf. \cite{HS1}). Also, based on
numerical evidence, in \cite{E1} Evans was led to propose a
conjecture which expresses $MK^7$ in terms of Hecke eigenvalues for
a weight 3 newform on $\Gamma_0(525)$ with quartic nebentypus of
conductor 105.

From now on, let us assume that $q=3^r$. Recently, Moisio was able
to find explicit expressions of $MK^h$, for $h\leq10$(cf.\cite{M1}).
This was done, via Pless power moment identity, by connecting
moments of Kloosterman sums and the frequencies of weights in the
ternary Melas code of length $q-1$, which were known by the work of
Geer, Schoof and Vlugt in \cite{GS1}.

 In order to describe our results, we introduce three incomplete
power moments of Kloosterman sums. For every nonnegative integer
$h$, and $\psi$ as before, we define

\begin{equation}\label{a}
T_{0}SK(\psi)^{h}=\sum_{a
\in\mathbb{F}_{q}^{*},tra=0}K(\psi;a^{2})^{h},~
T_{12}SK(\psi)^{h}=\sum_{a \in\mathbb{F}_{q}^{*},tra\neq
0}K(\psi;a^{2})^{h},
\end{equation}
which will be respectively  called the $h$-th moment of Kloosterman
sums with $\lq\lq$trace zero square arguments" and those with
$\lq\lq$trace nonzero square arguments." Then, clearly we have
\begin{equation}\label{b}
2SK(\psi)^{h}=T_{0}SK(\psi)^{h}+T_{12}SK(\psi)^{h},
\end{equation}
where
\begin{equation}\label{c}
SK(\psi)^{h}=\sum_{a \in\mathbb{F}_{q}^{*},a ~square
}K(\psi;a^{})^{h},
\end{equation}
called the $h$-th moment of Kloosterman sums with $\lq\lq$square
arguments." If $\psi=\lambda$ is the canonical additive character of
$\mathbb F_{q}$ ,then $SK(\lambda)^{h}$, $T_{0}SK(\lambda)^{h}$, and
$T_{12}SK(\lambda)^{h}$ will be respectively denoted by $SK^{h}$,
$T_{0}SK^{h}$, and $T_{12}SK^{h}$, for brevity.

We derived in \cite{D5} recursive formulas for the power moments of
Kloosterman sums with trace nonzero square arguments. To do that, we
constructed ternary linear codes $C(SO(3,q))$ and $C(O(3,q))$,
respectively associated with the
  orthogonal groups  $SO(3,q)$ and $O(3,q)$, and expressed those power moments in terms of the frequencies of weights in the codes.

 In this paper, we will obtain two infinite families of recursive
formulas, the one generating the power moments of Kloosterman sums
with trace nonzero square arguments and the other generating the
even power moments of those. To do that, we construct four infinite
families  of  ternary linear codes  associated with double cosets in
$O(2n+1,q)$  with respect to certain maximal parabolic subgroup of
the special orthogonal  group $SO(2n+1,q)$, and express those power
moments in terms of the frequencies of weights in the codes. Then,
thanks to our previous results  on the explicit expressions of
exponential sums over those double cosets related to the evaluations
of $\lq\lq$Gauss sums" for the orthogonal  group $O(2n+1,q)$
\cite{D3}, we can express the weight of each codeword in the duals
of the codes in terms of Kloosterman sums. Then our formulas will
follow immediately from the Pless power moment identity.

 The following Theorem \ref{A} is the main result of this paper.
Henceforth, we agree that, for nonnegative integers $a,b,c,$

\begin{equation}\label{d}
\binom{c}{a,b}=\frac{c!}{a!b!(c-a-b)!}, ~\text {if}~ a+b\leq c,
\end{equation}

\begin{equation}\label{e}
\binom{c}{a,b}=0, ~\text {if}~ a+b>c.
\end{equation}

 To simplify notations, we introduce the following ones which will be
used throughout this paper at various places.

\begin{equation}\label{f}
A^{-}(n,q)=q^{\frac{1}{4}(5n^{2}-1)}\Big[
                                      \begin{array}{c}
                                        n \\
                                        1 \\
                                      \end{array}
                                    \Big]
_{q}\prod^{(n-1)/2}_{j=1}(q^{2j-1}-1),
\end{equation}
\begin{equation}\label{g}
B^{-}(n,q)=q^{\frac{1}{4}(n-1)^{2}}(q^{n}-1)\prod^{(n-1)/2}_{j=1}(q^{2j}-1),\quad
\end{equation}
\begin{equation}\label{h}
A^{+}(n,q)=q^{\frac{1}{4}(5n^{2}-2n)}\Big[
                                       \begin{array}{c}
                                         n \\
                                         2 \\
                                       \end{array}
                                     \Big]
_{q}\prod^{(n-2)/2}_{j=1}(q^{2j-1}-1),
\end{equation}
\begin{equation}\label{i}
\quad\quad\quad
B^{+}(n,q)=q^{\frac{1}{4}(n-2)^{2}}(q^{n}-1)(q^{n-1}-1)\prod^{(n-2)/2}_{j=1}(q^{2j}-1).
\end{equation}

 From now on, it is assumed that either $+$ signs or $-$ signs are chosen
everywhere, whenever $\pm$ signs appear.

\begin{theorem}\label{A}
$(1)$ For each odd integer $n\geq1$ , all q , and h=1,2,3,$\cdots$,

\begin{equation}\label{j}
\begin{split}
&((-1)^{h+1}+2^{-h})T_{12}SK^{h}\\
&=-\sum^{h-1}_{j=1}((-1)^{j+1}+2^{-j})\binom{h}{j}B
^{-}(n,q)^{h-j}T_{12}SK^{j}\\
&\quad\quad\quad+qA^{-}(n,q)^{-h}\sum^{\min\{N_{i}^{-}(n,q),h\}}_{j=0}(-1)^{j}(C^{-}_{i,j}(n,q)-C^{-}_{j}(n,q))\\
&\quad\quad\quad\quad\times\sum
^{h}_{t=j}t!S(h,t)3^{h-t}2^{t-h-j}\binom{N^{-}_{i}(n,q)-j}{N^{-}_{i}(n,q)-t}(i=1,2),
\end{split}
\end{equation}
where $N_{i}^{-}(n,q)=|DC^{-}_{i}(n,q)|=A^{-}(n,q)B^{-}(n,q),$ for
$i=1,2,$ and $\{C^{-}_{1,j}(n,q)\}^{N^{-}_{1}(n,q)}_{j=0}$,
$\{C^{-}_{2,j}(n,q)\}^{N^{-}_{2}(n,q)}_{j=0},$
 and $\{C^{-}_{j}(n,q)\}^{N^{-}_{i}(n,q)}_{j=0}$  are respectively the
weight distributions of the ternary linear codes
$C(DC^{-}_{1}(n,q)),
 C(DC^{-}_{2}(n,q)),$ and $ C(DC^{-}(n,q))$ given by:

\begin{equation*}
C^{-}_{1,j}(n,q)=\sum\binom{q^{-1}A^{-}(n,q)(B^{-}(n,q)+1)}{\nu_{0},\mu_{0}}\binom{q^{-1}A^{-}(n,q)(B^{-}(n,q)+1)}{\nu_{2},\mu_{2}}
\end{equation*}
\begin{equation}\label{k}
\begin{split}
&\times \prod_{\beta^{2}-2\beta\neq0
~square}\binom{q^{-1}A^{-}(n,q)(B^{-}(n,q)+q+1)}{\nu_{\beta},\mu_{\beta}}\\
&\times \prod_{\beta^{2}-2\beta
~nonsquare}\binom{q^{-1}A^{-}(n,q)(B^{-}(n,q)-q+1)}{\nu_{\beta},\mu_{\beta}},
\end{split}
\end{equation}

\begin{equation*}
C^{-}_{2,j}(n,q)=\sum\binom{q^{-1}A^{-}(n,q)(B^{-}(n,q)+1)}{\nu_{0},\mu_{0}}\binom{q^{-1}A^{-}(n,q)(B^{-}(n,q)+1)}{\nu_{1},\mu_{1}}\\
\end{equation*}
\begin{equation}\label{l}
\begin{split}
&\times\prod_{\beta^{2}+2\beta\neq0
~square}\binom{q^{-1}A^{-}(n,q)(B^{-}(n,q)+q+1)}{\nu_{\beta},\mu_{\beta}}\\
&\times\prod_{\beta^{2}+2\beta
~nonsquare}\binom{q^{-1}A^{-}(n,q)(B^{-}(n,q)-q+1)}{\nu_{\beta},\mu_{\beta}},
\end{split}
\end{equation}

\begin{equation*}
C^{-}_{j}(n,q)=\sum\binom{q^{-1}A^{-}(n,q)(B^{-}(n,q)+1)}{\nu_{1},\mu_{1}}\binom{q^{-1}A^{-}(n,q)(B^{-}(n,q)+1)}{\nu_{-1},\mu_{-1}}\\
\end{equation*}
\begin{equation}\label{m}
\begin{split}
&\times\prod_{\beta^{2}-1\neq0
~square}\binom{q^{-1}A^{-}(n,q)(B^{-}(n,q)+q+1)}{\nu_{\beta},\mu_{\beta}}\\
&\times\prod_{\beta^{2}-1
~nonsquare}\binom{q^{-1}A^{-}(n,q)(B^{-}(n,q)-q+1)}{\nu_{\beta},\mu_{\beta}}(cf.
(\ref{d}),(\ref{e})) .
\end{split}
\end{equation}
Here the first sum in (\ref{j}) is $0$ if $h=1$ and the unspecified
sums in (\ref{k}), (\ref{l}), and (\ref{m}) run over all the sets of
nonnegative integers $\{\nu_{\beta}\}_{\beta\in\mathbb{F}_{q}}$, and
$\{\mu_{\beta}\}_{\beta\in\mathbb{F}_{q}}$  satisfying

\begin{align*}
\sum_{\beta\in\mathbb{F}_q} \nu_\beta
+\sum_{\beta\in\mathbb{F}_q}\mu_\beta=j, \quad \text{and}\quad
\sum_{\beta\in\mathbb{F}_q} \nu_\beta
\beta=\sum_{\beta\in\mathbb{F}_q}\mu_\beta \beta.
\end{align*}
In addition, $S(h,t)$  is the Stirling number of the second kind
defined by

\begin{equation}\label{n}
S(h,t)=\frac{1}{t!}\sum^{t}_{j=0}(-1)^{t-j}\binom{t}{j}j^{h}.
\end{equation}

\item\label{Ab}
$(2)$ For each even integer $n\geq2$, all $q$, and $h=1,2,3,\cdots,$
\begin{equation}\label{o}
\begin{split}
&((-1)^{h+1}+2^{-h})T_{12}SK^{2h}\\
&=-\sum^{h-1}_{j=0}\binom{h}{j}\{(-1)^{j+1}(B
^{+}(n,q)-q^{2}+q)^{h-j}\\
&\qquad\qquad+2^{-j}(B^{+}(n,q)+ \frac {1}{2}q^{2}-\frac{1}{2}q)^{h-j})\}T_{12}SK^{2j}\\
&\quad\quad+qA^{+}(n,q)^{-h}\sum^{\min\{N_{i}^{+}(n,q),h\}}_{j=0}(-1)^{j}(C^{+}_{i,j}(n,q)-C^{+}_{j}(n,q))\\
&\quad\quad\quad\quad\quad\times\sum
^{h}_{t=j}t!S(h,t)3^{h-t}2^{t-h-j}\binom{N^{+}_{i}(n,q)-j}{N^{+}_{i}(n,q)-t}~(i=1,2),
\end{split}
\end{equation}
where $N_{i}^{+}(n,q)=|DC^{+}_{i}(n,q)|=A^{+}(n,q)B^{+}(n,q), $ for
$i=1,2,$ and \\
$ \{C^{+}_{1,j}(n,q)\}^{N^{+}_{1}(n,q)}_{j=0}$,
$\{C^{+}_{2,j}(n,q)\}^{N^{+}_{2}(n,q)}_{j=0}, $ and
$\{C^{+}_{j}(n,q)\}^{N^{+}_{i}(n,q)}_{j=0}$  are respectively the
weight distributions of the ternary linear codes
$C(DC^{+}_{1}(n,q)),
 C(DC^{+}_{2}(n,q)),$ and $ C(DC^{+}(n,q))$ given by:
\begin{equation*}
\begin{split}
C^{+}_{1,j}(n,q)=&\sum\binom{q^{-1}A^{+}(n,q)(B^{+}(n,q)+q\delta(2,q;0)+(q-1)^{3})}{\nu_{1},\mu_{1}}\\
\end{split}
\end{equation*}
\begin{equation}\label{p}
\begin{split}
&\times \prod_{\beta\neq1
}\binom{q^{-1}A^{+}(n,q)(B^{+}(n,q)+q\delta(2,q;\beta-1)-2q^{2}+3q-1)}{\nu_{\beta},\mu_{\beta}},\\
\end{split}
\end{equation}
\begin{equation*}
\begin{split}
C^{+}_{2,j}(n,q)=&\sum \binom{q^{-1}A^{+}(n,q)(B^{+}(n,q)+q\delta(2,q;0)+(q-1)^{3})}{\nu_{-1},\mu_{-1}}\\
\end{split}
\end{equation*}
\begin{equation}\label{q}
\begin{split}
&\times \prod_{\beta\neq{-1}
}\binom{q^{-1}A^{+}(n,q)(B^{+}(n,q)+q\delta(2,q;\beta+1)-2q^{2}+3q-1)}{\nu_{\beta},\mu_{\beta}},\\
\end{split}
\end{equation}
\begin{equation*}
\begin{split}
C^{+}_{j}(n,q)=&\sum\binom{q^{4}(\delta(2,q;0)+q^{5}-q^{2}-3q+3)}{\nu_{0},\mu_{0}}\qquad\qquad\qquad\\
\end{split}
\end{equation*}
\begin{equation}\label{r}
\begin{split}
&\times \prod_{\beta\in \mathbb
{F}_{q}^{*}}\binom{q^{4}(\delta(2,q;\beta)+q^{5}-q^{3}-q^{2}-2q+3)}{\nu_{\beta},\mu_{\beta}}~(cf.(\ref{d}),(\ref{e})).\\
\end{split}
\end{equation}
Here the sums are over all the sets of nonnegative integers
$\{\nu_{\beta}\}_{\beta\in\mathbb{F}_{q}}$ and
$\{\mu_{\beta}\}_{\beta\in\mathbb{F}_{q}}$  satisfying

\begin{align*}
\sum_{\beta\in\mathbb{F}_q} \nu_\beta
+\sum_{\beta\in\mathbb{F}_q}\mu_\beta=j,\quad \text{and} \quad
\sum_{\beta\in\mathbb{F}_q} \nu_\beta
\beta=\sum_{\beta\in\mathbb{F}_q}\mu_\beta \beta,
\end{align*}
and $\delta(2,q;\beta)=|\{(\alpha_{1},\alpha_{2})\in
\mathbb{F}_{q}^{2}|\alpha_{1}+\alpha_{2}+\alpha_{1}^{-1}+\alpha_{2}^{-1}=\beta\}|$.
\end{theorem}

\section{$O(2n+1,q)$}
For more details about this section, one is referred to the paper
\cite{D3}. Throughout this paper, the following notations will be
used:
\begin{align*}
\begin{split}
q&=3^r~(r\in\mathbb{Z}_{>0}),\\
\mathbb{F}_q&=~the~finite~field~ with~ q~ elements,\\
TrA&=~the~ trace~ of~ A~ for~ a~ square~ matrix~ A,\\
^t B&=~the~transpose~ of~ B~for~any~matrix~B.
\end{split}
\end{align*}

The orthogonal group $O(2n+1,q)$  is defined as:
\begin{align*}
O(2n+1,q)=\{w\in GL(2n+1,q)|^twJw=J\},
\end{align*}
where
\begin{align*}
J=\begin{bmatrix}
    0 & 1_{n} & 0 \\
    1_{n} & 0 & 0 \\
    0 & 0 & 1 \\
  \end{bmatrix}
.
\end{align*}
It consists of the matrices

\begin{align*}
\begin{bmatrix}
 A & B & e \\
 C & D & f \\
 g & h & i \\
\end{bmatrix}
(A,B,C,D~n\times n,e,f~ n\times1,g,h~1\times n, i~ 1\times 1)
\end{align*}
in $GL(2n+1,q)$ satisfying the relations:

\begin{align*}
\begin{split}
&{}^tAC+{}^tCA+{}^tgg=0,\quad {}^tBD+{}^tDB+{}^thh=0,\\
&{}^tAD+{}^tCB+{}^tgh=1_{n},\quad {}^tef+{}^tfe+i^{2}=1,\\
&{}^tAf+{}^tCe+{}^tgi=0,\quad {}^tBf+{}^tDe+{}^thi=0.
\end{split}
\end{align*}

Let $P(2n+1,q)$ be the maximal parabolic subgroup of $O(2n+1,q)$
given by
\begin{align*}
\begin{split}
P&=P(2n+1,q)\\
&=\left\{\begin{bmatrix}
       A & 0 & 0 \\
        0 & ^tA^{-1} & 0 \\
        0 & 0 & i \\
    \end{bmatrix}
\begin{bmatrix}
 1_{n} & B & -^th \\
               0 & 1_{n} & 0 \\
               0 & h & 1 \\
\end{bmatrix}
\Bigg|\begin{array}{c}
   A\in GL(n,q), ~i=\pm 1 \\
   B+^tB+^thh=0
 \end{array}
\right\},
\end{split}
\end{align*}
and let $Q=Q(2n+1,q)$ be the subgroup of $P(2n+1,q)$ of index 2
defined by

\begin{equation*}
Q=Q(2n+1,q)\\
=\left\{\begin{bmatrix}
      A & 0 & 0 \\
      0 & ^tA^{-1} & 0 \\
      0 & 0 & 1 \\
    \end{bmatrix}
\begin{bmatrix}
 1_{n} & B & -^th \\
 0 & 1_{n} & 0 \\
 0 & h & 1 \\
\end{bmatrix}
\Bigg|\begin{array}{c}
   A\in GL(n,q) \\
   B+^tB+^thh=0
 \end{array}
\right\}.
\end{equation*}
Then we see that
\begin{equation}
P(2n+1,q)=Q(2n+1,q)\amalg\rho Q(2n+1,q),
\end{equation}
with
\begin{align*}
\rho=\begin{bmatrix}
        1_{n} & 0 & 0 \\
        0 & 1_{n} & 0 \\
        0 & 0 & -1 \\
     \end{bmatrix}
 .
\end{align*}

Let $\sigma_r$ denote the following matrix in $O(2n+1,q)$
\begin{align*}
\sigma_r=
\begin{bmatrix}
   0 & 0 & 1_r & 0 & 0 \\
   0 & 1_{n-r} & 0 & 0 & 0 \\
   1_r & 0 & 0 & 0 & 0 \\
   0 & 0 & 0 & 1_{n-r} & 0 \\
   0 & 0 & 0 & 0 & 1 \\
\end{bmatrix}
~(0\leq r\leq n).
\end{align*}
Then the Bruhat decomposition of $O(2n+1,q)$ with respect to
$P=P(2n+1,q)$ is given by
\begin{equation}\label{s}
\begin{split}
O(2n+1,q)&=\coprod_{r=0}^{n}P\sigma_r P=\coprod_{r=0}^{n}P\sigma_rQ\\
&=\coprod_{r=0}^{n}Q\sigma_r Q\amalg\coprod_{r=0}^{n}\rho
Q\sigma_rQ,
\end{split}
\end{equation}
which can further be modified as

\begin{equation}\label{t}
\begin{split}
O(2n+1,q)&=\coprod_{r=0}^{n}P\sigma_r (B_{r}\setminus Q)\\
&=\coprod_{r=0}^{n}Q\sigma_r (B_{r}\setminus
Q)\amalg\coprod_{r=0}^{n}\rho Q\sigma_r(B_{r}\setminus Q),
\end{split}
\end{equation}
with

\begin{align*}
B_r=B_r (q)=\{w\in Q(2n+1,q)|\sigma_r w \sigma_r^{-1}\in
P(2n+1,q)\}.
\end{align*}

  For integers $n$, $r$ with $0\leq r\leq n$, the $q$-binomial
coefficients are defined as:
\begin{align*}
\Big[
  \begin{array}{c}
    n \\
    r \\
  \end{array}
\Big] _q=\prod_{j=0}^{r-1}(q^{n-j}-1)/(q^{r-j}-1).
\end{align*}
It is shown in \cite{D3} that
\begin{equation}\label{u}
|B_{r}(q)\setminus Q(2n+1,q)|=q^{\binom {r+1}{2}} {\Big[
                                                     \begin{array}{c}
                                                       n \\
                                                       r \\
                                                     \end{array}
                                                   \Big]
 }_{q},
\end{equation}

\begin{equation}\label{v}
\begin{split}
|Q(2n+1,q)\sigma_{r}Q(2n+1,q)|&=|\rho
Q(2n+1,q)\sigma_{r}Q(2n+1,q)|\\
&=q^{n^{2}}\prod_{j=1}^{n}(q^{j}-1)q^{\binom {r}{2}}q^{r}{\Big[
                                                            \begin{array}{c}
                                                              n \\
                                                              r \\
                                                            \end{array}
                                                          \Big]
 }_{q}.
\end{split}
\end{equation}

Let
\begin{equation}\label{w}
DC_{1}^{-}(n,q)=Q(2n+1,q)\sigma_{n-1}Q(2n+1,q), \quad \text { ~for~
} n=1,3,5,\cdots,
\end{equation}

\begin{equation}\label{x}
DC_{2}^{-}(n,q)=\rho Q(2n+1,q)\sigma_{n-1}Q(2n+1,q), \quad \text {
~for~ }~ n=1,3,5,\cdots,
\end{equation}

\begin{equation}\label{y}
DC_{1}^{+}(n,q)=Q(2n+1,q)\sigma_{n-2}Q(2n+1,q), \quad \text { ~for~
}~ n=2,4,6,\cdots,
\end{equation}

\begin{equation}\label{z}
DC_{2}^{+}(n,q)=\rho Q(2n+1,q)\sigma_{n-2}Q(2n+1,q), \quad \text {
~for~ }~ n=2,4,6,\cdots.
\end{equation}
Then, from (\ref {v}), we have:

\begin{equation}\label{a1}
|DC_{i}^{\mp}(n,q)|=A^{\mp}(n,q)B^{\mp}(n,q), \quad \text {for}~
i=1,2~\text {(cf. $($\ref{f}$)$-$($\ref{i}$)$)}.
\end{equation}
\\
 Unless otherwise stated, from now on, we will agree that anything
related to $DC_{1}^{-}(n,q)$ and  $DC_{2}^{-}(n,q)$ are defined for
$n=1,3,5,\cdots ,$ and anything related to  $DC_{1}^{+}(n,q)$ and
$DC_{2}^{+}(n,q)$ are defined for $n=2,4,6,\cdots .$


\section{Exponential sums over double cosets of  $O(2n+1,q)$}

 The following notations will be employed throughout this paper.
\begin{align*}
\begin{split}
tr(x)&=x+x^3+\cdots+x^{3^{r-1}}~the~trace~function~\mathbb{F}_q\rightarrow\mathbb{F}_3,\\
 \lambda_0(x)&=e^{2\pi
ix/3}~the~canonical~additive~character~of~\mathbb{F}_3,\\
\lambda(x)&=e^{2\pi
tr(x)/3}~the~canonical~additive~character~of~\mathbb{F}_q.
\end{split}
\end{align*}
Then any nontrivial additive character $\psi$ of $\mathbb{F}_q$ is
given by $\psi(x)=\lambda(ax)$, for a unique
$a\in\mathbb{F}_{q}^{*}$. Also, since $\lambda(a)$ for any $a\in
\mathbb{F}_q$ is a 3th root of 1, we have

\begin{equation}\label{b1}
\lambda(-a)=\lambda(2a)=\lambda(a)^{2}=\lambda(a)^{-1}=\overline{\lambda(a)}.
\end{equation}

For any nontrivial additive character $\psi$ of $\mathbb{F}_q$ and
$a\in\mathbb{F}_{q}^{*}$, the Kloosterman sum $K_{GL(t,q)}(\psi ;
a)$ for $GL(t,q)$ is defined as

\begin{align*}
K_{GL(t,q)}(\psi;a)=\sum_{w\in GL(t,q)}\psi(Trw+aTrw^{-1}).
\end{align*}
Observe  that, for $t=1$ , $K_{GL(1,q)}(\psi;a)$ denotes the
Kloosterman sum $K(\psi;a)$.
In \cite{D1}, it is shown that
$K_{GL(t,q)}(\psi;a)$ satisfies the following recursive relation:
for integers $t\geq2$, $a\in \mathbb{F}_{q}^{*}$,

\begin{align*}
K_{GL(t,q)}(\psi;a)=q^{t-1}K_{GL(t-1,q)}(\psi;a)K(\psi;a)+q^{2t-2}(q^{t-1}-1)K_{GL(t-2,q)}(\psi;a),
\end{align*}
where we understand that $K_{GL(0,q)}(\psi,a)=1$.

\begin{proposition}\label{B}$($\cite{D3}$)$
Let $\psi$ be a nontrivial additive character of $\mathbb{F}_q$. For
each positive integer $r$, let $\Omega_r$ be the set of all $r\times
r$ nonsingular symmetric matrices over $\mathbb F_q$. Then we have

\begin{equation}\label{c1}
a_r(\psi)=\sum_{B\in\Omega_r}\sum_{h\in\mathbb{F}_{q}^{r\times
1}}\psi({^t}hBh) =
\begin{cases}
\begin{split}
q^{r(r+2)/4}\prod_{j=1}^{r/2}(q^{2j-1}-1),& \quad \text {for r even,}\\
0,& \quad \text {for r odd.}
\end{split}
\end{cases}
\end{equation}
\end{proposition}
From Sections 5 and 6 of \cite{D3} the Gauss sum for $O(2n+1,q)$,
with $\psi$ a nontrivial additive character of $\mathbb F_q$, is
given by:
\begin{align*}
\sum_{w\in O(2n+1,q)}\psi(Trw) =\sum_{0\leq r \leq n}\sum_{w\in
Q\sigma_{r}Q}\psi(Trw) +\sum_{0\leq r \leq n}\sum_{w\in \rho
Q\sigma_{r}Q}\psi(Trw)~(cf.~ (\ref{s})),
\end{align*}
with
\begin{equation}\label{d1}
\begin{split}
&\sum_{w\in Q\sigma_{r}Q}\psi(Trw)\\
&=|B_r\setminus Q|\sum_{w\in Q}\psi(Trw\sigma_{r})\\
&=\psi(1)q^{\binom{n+1}{2}}|B_r\setminus
Q|q^{r(n-r-1)}a_{r}(\psi)K_{GL(n-r,q)}(\psi;1),
\end{split}
\end{equation}
\begin{equation}\label{e1}
\begin{split}
&\sum_{w\in \rho Q\sigma_{r}Q}\psi(Trw)\\
&=|B_r\setminus Q|\sum_{w\in Q}\psi(Tr \rho w\sigma_{r})\\
&=\psi(-1)q^{\binom{n+1}{2}}|B_r\setminus
Q|q^{r(n-r-1)}a_{r}(\psi)K_{GL(n-r,q)}(\psi;1).
\end{split}
\end{equation}
Here one uses (\ref{t}) and the fact that $\rho^{-1}w\rho \in Q$,
for all $w\in Q$.\\

We now see from (\ref{u}) and (\ref{c1})-(\ref{e1}) that, for each
$r$ with $0\leq r\leq n$,

\begin{equation}\label{f1}
\sum_{w\in
Q\sigma_{r}Q}\psi(Trw)\qquad\qquad\qquad\qquad\qquad\qquad\qquad\qquad\qquad\qquad\qquad\qquad
\end{equation}
\begin{equation*}
=\begin{cases}
\psi(1)q^{\binom{n+1}{2}}q^{rn-\frac{1}{4}r^{2}}{\Big[
                                                   \begin{array}{c}
                                                     n \\
                                                     r \\
                                                   \end{array}
                                                 \Big]
}_{q}\prod_{j=1}^{r/2}(q^{2j-1}-1)K_{GL(n-r,q)}(\psi;1),& ~if~ r ~is ~even,\\
0,& ~if ~r~is~ odd,
\end{cases}
\end{equation*}

\begin{equation}\label{g1}
\sum_{w\in \rho
Q\sigma_{r}Q}\psi(Trw)\qquad\qquad\qquad\qquad\qquad\qquad\qquad\qquad\qquad\qquad\qquad\qquad
\end{equation}
\begin{equation*}
=\begin{cases}
\psi(-1)q^{\binom{n+1}{2}}q^{rn-\frac{1}{4}r^{2}}{\Big[
                                                   \begin{array}{c}
                                                     n \\
                                                     r \\
                                                   \end{array}
                                                 \Big]}_{q}\prod_{j=1}^{r/2}(q^{2j-1}-1)K_{GL(n-r,q)}(\psi;1),& ~if~ r ~is ~even,\\
0,& ~if~ r~ is~ odd.
\end{cases}
\end{equation*}

For our purposes, we need two infinite families of exponential sums
in (\ref{f1}) over  $DC_{1}^{-}(n,q)$,   for $n=1,3,5,\cdots,$ and
over $DC_{1}^{+}(n,q)$, for $n=2,4,6,\cdots$. And also, we need two
such sums in (\ref{g1}) over $DC_{2}^{-}(n,q)$, for
$n=1,3,5,\cdots,$ and over $DC_{2}^{+}(n,q)$, for $n=2,4,6,\cdots$.
So we state them separately as a theorem.

\begin{theorem}\label{C}
Let $\psi$  be any nontrivial additive character of
$\mathbb{F}_{q}$. Then, in the notations of (\ref{f}), (\ref{h}),
and (\ref{w})-(\ref{z}), we have

\begin{equation*}
\begin{split}
\sum_{w\in DC_{1}^{-}(n,q)}\psi(Trw)&=\psi(1)A^{-}(n,q)K(\psi;1),
~\text{for}~ n=1,3,5,\cdots,\\
\sum_{w\in DC_{2}^{-}(n,q)}\psi(Trw)&=\psi(-1)A^{-}(n,q)K(\psi;1),
~\text{for}~ n=1,3,5,\cdots,\\
\sum_{w\in
DC_{1}^{+}(n,q)}\psi(Trw)&=\psi(1)q^{-1}A^{+}(n,q)K_{GL(2,q)}(\psi;1)\\
                         &=\psi(1)A^{+}(n,q)(K(\psi;1)^{2}+q^{2}-q),
                         ~\text{for}~ n=2,4,6,\cdots,\\
\sum_{w\in
DC_{2}^{+}(n,q)}\psi(Trw)&=\psi(-1)q^{-1}A^{+}(n,q)K_{GL(2,q)}(\psi;1)\\
                         &=\psi(-1)A^{+}(n,q)(K(\psi;1)^{2}+q^{2}-q),
                         ~\text{for}~ n=2,4,6,\cdots.
\end{split}
\end{equation*}
\end{theorem}

The next corollary follows from Theorem \ref{C}  and simple changes
of variables.

\begin{corollary}\label{D}
 Let $\lambda$ be the canonical additive character of $\mathbb{F}_{q}$,
and let $a\in \mathbb{F}_{q}^{*}$. Then we have
\begin{equation}\label{h1}
\sum_{w \in DC_{1}^{-}(n,q)}
\lambda(aTrw)=\lambda(a)A^{-}(n,q)K(\lambda;a^{2}), ~\text{for}~
n=1,3,5,\cdots,
\end{equation}
\begin{equation}\label{i1}
\sum_{w \in DC_{2}^{-}(n,q)}
\lambda(aTrw)=\lambda(-a)A^{-}(n,q)K(\lambda;a^{2}), ~\text{for}~
n=1,3,5,\cdots,
\end{equation}
\begin{equation}\label{j1}
\sum_{w \in DC_{1}^{+}(n,q)}
\lambda(aTrw)=\lambda(a)A^{+}(n,q)(K(\lambda;a^{2})^{2}+q^2-q),
~\text{for}~ n=2,4,6,\cdots,
\end{equation}
\begin{equation}\label{k1}
\sum_{w \in DC_{2}^{+}(n,q)}
\lambda(aTrw)=\lambda(-a)A^{+}(n,q)(K(\lambda;a^{2})^{2}+q^2-q),
~\text{for}~ n=2,4,6,\cdots.
\end{equation}
\end{corollary}

\begin{proposition}\label{E}$([5, (5.3$-$5)])$
Let $\lambda$ be the canonical additive character of $\mathbb{F}_q$,
$m\in \mathbb{Z}_{\geq0}$, $\beta\in\mathbb{F}_q$. Then
\begin{equation}\label{l1}
\sum_{a\in\mathbb{F}_{q}^{*}}\lambda(-a\beta)K(\lambda;a^2)^m=q\delta(m,q;\beta)-(q-1)^m,
\end{equation}
where, for $m\geq1$,
\begin{equation}\label{m1}
\delta(m,q;\beta)=|\{(\alpha_{1}^{},\cdots,\alpha_{m}^{})\in(\mathbb{F}_{q}^{*})^m|\alpha_{1}^{}+\alpha_{1}^{-1}+
\cdots,+\alpha_{m}^{}+\alpha_{m}^{-1}=\beta\}|,
\end{equation}
and
\begin{equation}\label{n1}
\delta(0,q;\beta)=
\begin{cases}
1,& ~if~\beta=0,\\
0,& ~otherwise.
\end{cases}
\end{equation}
\end{proposition}

\begin{remark}\label{F}
Here one notes that
\begin{equation}\label{o1}
\begin{split}
\delta(1,q;\beta)&=|\{x\in\mathbb{F}_{q}^{}|x^2-\beta x+1=0\}|\\
&=
\begin{cases}
2,& \text {if $\beta^2-1\neq$0 is a square,}\\
1,& \text {if $\beta^2-1$=0,}\\
0,& \text {if $\beta^2-1$ is a nonsquare.}
\end{cases}
\end{split}
\end{equation}
\end{remark}

\begin{lemma}\label{G}$($\cite{K1}$)$
Let $\delta(m,q;\beta)$ be as in  (\ref{m1}) and (\ref{n1}), and let
$a \in \mathbb{F}_{q}^{*}$. Then we have

\begin{equation}\label{p1}
\sum_{\beta \in
\mathbb{F}_{q}}\delta(m,q;\beta)\lambda(a\beta)=K(\lambda;a^{2})^{m}.
\end{equation}
 For any integer $r$ with $0\leq r\leq n$, and each $\beta \in
\mathbb{F}_{q}$, we let

\begin{equation*}
\begin{split}
N_{Q\sigma_{r}Q}(\beta)=|\{w\in Q\sigma_{r}Q|~Trw=\beta\}|,\\
N_{\rho Q\sigma_{r}Q}(\beta)=|\{w\in \rho
Q\sigma_{r}Q|~Trw=\beta\}|.
\end{split}
\end{equation*}
Then it is easy to see that
\begin{equation}\label{q1}
qN_{Q\sigma_{r}Q}(\beta)=|Q\sigma_{r}Q|+\sum_{a\in
\mathbb{F}_{q}^{*}}\lambda(-a\beta)\sum_{w\in
Q\sigma_{r}Q}\lambda(aTrw),
\end{equation}

\begin{equation}\label{r1}
qN_{\rho Q\sigma_{r}Q}(\beta)=|\rho Q\sigma_{r}Q|+\sum_{a\in
\mathbb{F}_{q}^{*}}\lambda(-a\beta)\sum_{w\in \rho
Q\sigma_{r}Q}\lambda(aTrw).
\end{equation}
\end{lemma}

Now, from (\ref{w})-(\ref{a1}) and (\ref{h1})-(\ref{k1}), we have
the following result.

\begin{proposition}\label{E}
With the notations of (\ref{f})-(\ref{i}), we have:
\item\label{Ha}
$(1)$
\begin{equation*}
\begin{split}
N_{DC_{1}^{-}(n,q)}&(\beta)=q^{-1}A^{-}(n,q)B^{-}(n,q)+q^{-1}A^{-}(n,q)(q\delta(1,q;\beta-1)-q+1)\\
&=q^{-1}A^{-}(n,q)B^{-}(n,q)+q^{-1}A^{-}(n,q)\\
\end{split}
\end{equation*}
\begin{equation}\label{s1}
\begin{split}
&\quad\quad\quad\quad\quad\times
\begin{cases}
q+1,& \text {~if~} \beta^{2}-2\beta\neq0 ~is~a~square~,\\
1,& \text {~if~} \beta=0 ~or~ 2,\\
-q+1, & \text {~if~} \beta^{2}-2\beta ~is ~a ~nonsquare,
\end{cases}
~(cf.(\ref{m1}),(\ref{o1}))
\end{split}
\end{equation}
\item\label{Hb}
$(2)$
\begin{equation*}
\begin{split}
N_{DC_{2}^{-}(n,q)}&(\beta)=q^{-1}A^{-}(n,q)B^{-}(n,q)+q^{-1}A^{-}(n,q)(q\delta(1,q;\beta+1)-q+1)\\
&=q^{-1}A^{-}(n,q)B^{-}(n,q)+q^{-1}A^{-}(n,q)\\
\end{split}
\end{equation*}
\begin{equation}\label{t1}
\begin{split}
&\times
\begin{cases}
q+1,& \text {~if~} \beta^{2}+2\beta\neq0 ~is ~a ~square,\\
1,& \text {~if~} \beta=0 ~or~ 1,\\
-q+1, & \text {~if~} \beta^{2}+2\beta  ~is ~a ~nonsquare,
\end{cases}
\end{split}
\end{equation}
\item\label{Hc}
$(3)$
\begin{equation}\label{u1}
\begin{split}
N_{DC_{1}^{+}(n,q)}&(\beta)=q^{-1}A^{+}(n,q)B^{+}(n,q)\\
&+q^{-1}A^{+}(n,q)\times
\begin{cases}
q\delta(2,q;0)+(q-1)^{3},& \text {~if~} \beta=1,\\
q\delta(2,q;\beta-1)-2q^{2}+3q-1,& \text {~if~} \beta\neq1,
\end{cases}
\end{split}
\end{equation}
(cf.(\ref{m1}))
\item\label{Hd}
$(4)$
\begin{equation}\label{v1}
\begin{split}
N_{DC_{2}^{+}(n,q)}&(\beta)=q^{-1}A^{+}(n,q)B^{+}(n,q)\\
&+q^{-1}A^{+}(n,q)\times
\begin{cases}
q\delta(2,q;0)+(q-1)^{3},& \text {~if~} \beta=-1,\\
q\delta(2,q;\beta+1)-2q^{2}+3q-1,& \text {~if~} \beta\neq-1.
\end{cases}
\end{split}
\end{equation}
\end{proposition}

\begin{corollary}\label{I}
$(1)$ For each odd integer $n\geq3$, with all $q$,
$N_{DC_{1}^{-}(n,q)}(\beta)>0$, for all $\beta$; for $n=1$, with all
$q$,
\begin{equation}\label{w1}
\begin{split}
N_{DC_{1}^{-}(n,q)}&(\beta)=
\begin{cases}
2q,& \text {~if~} \beta^{2}-2\beta\neq0 \text{~is~a~square~},\\
q,& \text {~if~} \beta=0 \text{~or~} 2,\\
0, & \text {~if~} \beta^{2}-2\beta \text{~is a nonsquare}.
\end{cases}
\end{split}
\end{equation}

\item\label{Ib}
$(2)$ For each odd integer $n\geq3$, with all $q$,
$N_{DC_{2}^{-}(n,q)}(\beta)>0$, for all $\beta$; for $n=1$, with all
$q$,
\begin{equation}\label{x1}
\begin{split}
N_{DC_{2}^{-}(n,q)}&(\beta)=
\begin{cases}
2q,& \text {~if~} \beta^{2}+2\beta\neq0 \text{~is~a~square~},\\
q,& \text {~if~} \beta=0 \text{~or~} 1,\\
0, & \text {~if~} \beta^{2}+2\beta \text{~is a nonsquare}.
\end{cases}
\end{split}
\end{equation}
\item\label{Ic}
$(3)$ For each even integer $n\geq2$, with all $q$,
$N_{DC_{1}^{+}(n,q)}(\beta)>0$, for all $\beta$.
\item\label{Id}
$(4)$ For each even integer $n\geq2$, with all $q$,
$N_{DC_{2}^{+}(n,q)}(\beta)>0$, for all $\beta$.
\end{corollary}

\section{Construction of codes}
Here we will construct two infinite families of ternary linear codes
$C(DC_{i}^{\mp}(n,q))$ of length $A^{\mp}(n,q)B^{\mp}(n,q)$,
associated with the double cosets $DC_{i}^{\mp}(n,q)~(i=1,2)$.

Let
\begin{equation}\label{y1}
N_{i}^{-}(n,q)=|DC_{i}^{-}(n,q)|=A^{-}(n,q)B^{-}(n,q),~\text{for}~i=1,2,
~\text{and}~ n=1,3,5,\cdots,
\end{equation}
\begin{equation}\label{z1}
N_{i}^{+}(n,q)=|DC_{i}^{+}(n,q)|=A^{+}(n,q)B^{+}(n,q),
~\text{for}~i=1,2,~\text{and}~ n=2,4,6,\cdots
\end{equation}
(cf.(\ref{a1})).

Let $g_{1},g_{2},\cdots,g_{N_{i}^{\mp}}(n,q)$ be fixed orderings of
the elements in $DC_{i}^{\mp}(n,q)$, for $ i=1,2,$ by abuse of
notations. Then we put
\begin{equation*}
v_{i}^{\mp}(n,q)=(Trg_1,Trg_2,\cdots,Trg_{N_{i}^{\mp}(n,q)}) \in
\mathbb F _{q}^{N_{i}^{\mp}(n,q)}, \quad \text{for}~i=1,2.
\end{equation*}
The ternary linear codes $C(DC_{1}^{-}(n,q))$ , $C(DC_{2}^{-}(n,q))$
, $C(DC_{1}^{+}(n,q))$ and\\
$C(DC_{2}^{+}(n,q))$  are defined as:

\begin{equation}\label{a2}
C(DC_{i}^{\mp}(n,q))=\{u\in \mathbb F
_{2}^{N_{i}^{\mp}(n,q)}|~u\cdot v_{i}^{\mp}(n,q)=0\}, ~\text{for}~
i=1,2,
\end{equation}
where the dot denotes respectively the usual inner product in
$\mathbb F_{q}^{N_{i}^{\mp}(n,q)}$, for $i=1,2.$

The following Delsarte's theorem is well-known.

\begin{theorem}\label{J}$($\cite{MS1}$)$
 Let $B$ be a linear code over $\mathbb{F}_{q}$. Then
\begin{equation*}
(B|_{\mathbb{F}_{3}})^{\perp}=tr(B^{\perp}).
\end{equation*}

  In view of this theorem, the dual  $C(DC_{i}^{\mp}(n,q))^{\perp}$ of the
code $C(DC_{i}^{\mp}(n,q))$ is given by, for $i=1,2,$

\begin{equation}\label{b2}
C(DC_{i}^{\mp}(n,q))^{\perp}=\{c_{i}^{\mp}(a)=c_{i}^{\mp}(a;n,q)=
(tr(aTrg_1),\cdots,tr(aTrg_{N_{i}^{\mp}(n,q)}))|a\in \mathbb{F}
_{q}\}.
\end{equation}
\end{theorem}
\begin{theorem}\label{K}
$(1)$ The map $\mathbb F_{q}\rightarrow
C(DC_{i}^{-}(n,q))^{\perp}(a\mapsto c_{i}^{-}(a))(i=1,2)$ is an
$\mathbb F_{3}$-linear isomorphism for each odd integer $n\geq1$ and
all $q$.
\item\label{Kb}
$(2)$ The map $\mathbb F_{q}\rightarrow
C(DC_{i}^{+}(n,q))^{\perp}(a\mapsto c_{i}^{+}(a))(i=1,2)$ is an
$\mathbb F_{3}$-linear isomorphism for each even integer $n\geq2$
and all $q$.
\begin{proof}
 All maps are clearly $\mathbb F_{3}$-linear and surjective. Let $a$ be
 in the kernel of map $\mathbb F_{q}\rightarrow
C(DC_{1}^{+}(n,q))^{\perp}(a\mapsto c_{1}^{+}(a))$. Then
$tr(aTrg)=0$, for all $g\in DC_{1}^{+}(n,q)$. Since, by Corollary
3.8 (3), $Tr:DC_{1}^{+}(n,q)\rightarrow \mathbb F_{q}$ is
surjective, and hence $tr(a\alpha)=0$, for all $\alpha \in\mathbb
F_{q}$. This implies that $a=0$, since otherwise $tr:\mathbb
F_{q}\rightarrow\mathbb F_{3}$ would be the zero map. This shows
$i=1$ case of  $(2)$. All the other assertions can be handled in the
same way, except for $i=1,2$ and $n=1$ case  of $(1)$, since in
those cases the maps $Tr:DC_{i}^{\mp}(n,q)\rightarrow \mathbb F_{q}$
are surjective.

Let $a$ be in the kernel of the map $\mathbb F_{q}\rightarrow
C(DC_{i}^{-}(1,q))^{\perp}(a\mapsto c_{i}^{-}(a))$, for $i=1,2.$
Then $tr(aTrg)=0$, for all $g\in DC_{i}^{-}(1,q)$. Suppose that
$a\neq0$. Then we would have\\
\begin{equation*}
\begin{split}
q(q-1)=|DC_{i}^{-}(1,q)|&=\sum_{g\in DC_{i}^{-}(1,q)}e^{2\pi itr(aTrg)/3}\\
&=\sum_{\beta \in \mathbb F_{q} }N_{DC_{i}^{-}(1,q)}(\beta)\lambda(a\beta)\\
&=q\sum_{\beta \in \mathbb F_{q}}\delta(1,q;\beta \mp1)\lambda(a\beta)(cf.(\ref{s1}),(\ref{t1}))\\
&(\text{Note here that it is $\beta-1$, for $i=1$, and $\beta+1$, for $i=2$})\\
&=q\lambda(\pm a)\sum_{\beta \in \mathbb F_{q}}\delta(1,q;\beta \mp1)\lambda(a(\beta\mp1))\\
&=q\lambda(\pm a)\sum_{\beta \in \mathbb F_{q}}\delta(1,q;\beta)\lambda(a\beta)\\
&=q\lambda(\pm a)K(\lambda;a^2)(cf.(\ref{p1})).\\
\end{split}
\end{equation*}
So, using Weil bound in (\ref{a}), we would get
\begin{equation*}
q-1=|K(\lambda;a^2)|\leq2\sqrt{q}.
\end{equation*}
For $q\geq9$, this is impossible. On the other hand, from
(\ref{w1})and (\ref{x1}) we see that, out of 6 elements in
$DC_{1}^{-}(1,3)$, 3 of them has $Tr=0$ and 3 of them has $Tr=2$;
out of 6 elements in $DC_{2}^{-}(1,3)$, 3 of them has $Tr=0$ and 3
of them has $Tr=1$. So in either case the kernel is trivial.
\end{proof}
\end{theorem}

\section{Power moments of Kloosterman sums with trace nonzero square arguments}
Here we will be able to find, via Pless power moment identity, two
infinite families of recursive formulas, the one generating the
power moments of  Kloosterman sums with trace nonzero square
arguments and the other generating the even power moments of those.

\begin{theorem}\label{L}$($Pless power moment identity,
\cite{MS1}$)$ Let $B$ be an $q$-ary $[$n,k$]$ code, and let $ B_{i}$
$($resp.$B_{i}^{\perp}$$)$ denote the number of codewords of weight
$i$ in $B$$($resp. in $B_{}^{\perp}$$)$.  Then, for
$h=0,1,2,\cdots,$
\begin{equation}\label{c2}
\sum_{j=0}^{n}j^{h}B_{j}=\sum_{j=0}^{min\{n,h\}}(-1)^{j}B_{j}^{\perp}\sum_{t=j}^{h}t!S(h,t)q^{k-t}(q-1)^{t-j}\binom{n-j}{n-t},
\end{equation}
where $S(h,t)$ is the Stirling number of the second kind defined in
(\ref{n}).
\end{theorem}

\begin{lemma}\label{M}
Let
$c_{i}^{\mp}(a)=(tr(aTrg_{1}),\cdots,tr(aTrg_{N_{i}^{\mp}(n,q)}))\in
C(DC_{i}^{\mp}(n,q))^{\perp}$, for $a \in \mathbb{F}_{q}^{*}$, and
 $i=1,2.$ Then the Hamming weights $w(c_{i}^{\mp}(a))$ are expressed
as follows:
\item\label{Ma}
$(1)$ $\qquad$ $w(c_{i}^{-}(a))$
\begin{equation}\label{d2}
=\frac{2}{3}A^{-}(n,q)(B^{-}(n,q)-(Re\lambda(a))K(\lambda;a^{2})),
~for~ i=1,2,
\end{equation}
\item\label{Mb}
$(2)$ $\qquad$ $w(c_{i}^{+}(a))$
\begin{equation}\label{e2}
\quad=\frac{2}{3}A^{+}(n,q)\{B^{+}(n,q)-(Re\lambda(a))(K(\lambda;a^{2})^{2}+q^2-q)\},
~for~ i=1,2.
\end{equation}
\begin{proof}
\begin{equation}
\begin{split}
w(c_{i}^{\mp}(a))&=\sum_{j=1}^{N_{i}^{\mp}(n,q)}(1-\frac{1}{3}\sum_{\alpha
\in \mathbb{F}_{3}^{}}\lambda_{0}(\alpha tr(aTrg_{j})))\\
&=N_{i}^{\mp}(n,q)-\frac{1}{3}\sum_{\alpha\in \mathbb
F_{3}}\sum_{w\in DC_{i}^{\mp}(n,q)}\lambda(\alpha aTrw)\\
&=\frac{2}{3}N_{i}^{\mp}(n,q)-\frac{1}{3}\sum_{\alpha\in \mathbb
F_{3}^{*}}\sum_{w\in DC_{i}^{\mp}(n,q)}\lambda(\alpha aTrw).
\end{split}
\end{equation}
Now, the results follow from (\ref{h1}), (\ref{k1}), (\ref{y1}), and
(\ref{z1}).
\end{proof}
\end{lemma}
 Let $u=(u_{1},\cdots,u_{N_{i}^{\mp}(n,q)})\in
\mathbb F_{3}^{N_{i}^{\mp}(n,q)}$, with $\nu_{\beta}$ 1's and
$\mu_{\beta}$ 2's in the coordinate places where $Tr(g_{j})=\beta$,
for each $\beta\in\mathbb F_{q}$. Then we see from the definition of
the code $C(DC_{i}^{\mp}(n,q))$(cf. (\ref{a2})) that $u$ is a
codeword with weight $j$ if and only if
$\sum_{\beta\in\mathbb{F}_q}\nu_\beta+\sum_{\beta\in\mathbb{F}_q}\mu_\beta=j$
and $\sum_{\beta\in\mathbb{F}_q}\nu_\beta
\beta=\sum_{\beta\in\mathbb{F}_q}\mu_\beta \beta$ (an identity in
$\mathbb{F}_q$). Note that there are
$\prod_{\beta\in\mathbb{F}_q}{\binom{N_{DC_{i}^{\mp}(n,q)}(\beta)}{\nu_\beta,\mu_\beta}}$(cf.
(\ref{d}), (\ref{e})) many such codewords with weight $j$. Now, we
get the following formulas in (\ref{f2})-(\ref{i2}), by using the
explicit values of $N_{DC_{i}^{\mp}(n,q)}(\beta)$ in
(\ref{s1})-(\ref{v1})(cf. (\ref{f})-(\ref{i})).

\begin{theorem}\label{N}
Let $\{C_{i,j}^{\mp}(n,q)\}_{j=0}^{N_{i}^{\mp}(n,q)}$  be the weight
distribution of  $C(DC_{i}^{\mp}(n,q))$, for $i=1,2.$ Then we
have:\\
\item\label{Na}
$(1)$ For $j=0,\cdots,N_{1}^{-}(n,q)$,
\begin{equation*}
\begin{split}
C_{1,j}^{-}(n,q)&=\sum_{\beta\in\mathbb{F}_q}\prod_{\beta\in
\mathbb{F}_q}
\binom{q^{-1}A^{-}(n,q)(B^{-}(n,q)+q\delta(1,q;\beta-1)-q+1)}{\nu_{\beta},\mu_{\beta}}\\
\end{split}
\end{equation*}
\begin{equation}\label{f2}
\begin{split}
&=\sum\binom{q^{-1}A^{-}(n,q)(B^{-}(n,q)+1)}{\nu_{0},\mu_{0}}
\binom{q^{-1}A^{-}(n,q)(B^{-}(n,q)+1)}{\nu_{2},\mu_{2}}\\
&\quad\times \prod_{\beta^{2}-2\beta\neq
0~square~}\binom{q^{-1}A^{-}(n,q)(B^{-}(n,q)+q+1)}{\nu_{\beta},\mu_{\beta}}\\
&\quad\times \prod_{\beta^{2}-2\beta
~nonsquare~}\binom{q^{-1}A^{-}(n,q)(B^{-}(n,q)-q+1)}{\nu_{\beta},\mu_{\beta}},
\end{split}
\end{equation}\\
\item\label{Nb}
$(2)$ For $j=0,\cdots,N_{2}^{-}(n,q)$,
\begin{equation*}
\begin{split}
C_{2,j}^{-}(n,q)&=\sum_{\beta\in\mathbb{F}_q}\prod_{\beta\in
\mathbb{F}_q}
\binom{q^{-1}A^{-}(n,q)(B^{-}(n,q)+q\delta(1,q;\beta+1)-q+1)}{\nu_{\beta},\mu_{\beta}}\\
\end{split}
\end{equation*}
\begin{equation}\label{g2}
\begin{split}
&=\sum\binom{q^{-1}A^{-}(n,q)(B^{-}(n,q)+1)}{\nu_{0},\mu_{0}}
\binom{q^{-1}A^{-}(n,q)(B^{-}(n,q)+1)}{\nu_{1},\mu_{1}}\\
&\quad\times \prod_{\beta^{2}+2\beta\neq
0~square~}\binom{q^{-1}A^{-}(n,q)(B^{-}(n,q)+q+1)}{\nu_{\beta},\mu_{\beta}}\\
&\quad\times \prod_{\beta^{2}+2\beta
~nonsquare~}\binom{q^{-1}A^{-}(n,q)(B^{-}(n,q)-q+1)}{\nu_{\beta},\mu_{\beta}}
(cf.(\ref{m1}),(\ref{o1})),
\end{split}
\end{equation}\\
\item\label{Nc}
$(3)$ For $j=0,\cdots,N_{1}^{+}(n,q)$,
\begin{equation*}
\begin{split}
C_{1,j}^{+}(n,q)&=\sum
\binom{q^{-1}A^{+}(n,q)(B^{+}(n,q)+q\delta(2,q;0)+(q-1)^3)}{\nu_{1},\mu_{1}}\\
\end{split}
\end{equation*}
\begin{equation}\label{h2}
\begin{split}
&\quad\times \prod_{\beta\neq
1}\binom{q^{-1}A^{+}(n,q)(B^{+}(n,q)+q\delta(2,q;\beta-1)-2q^2+3q-1)}{\nu_{\beta},\mu_{\beta}},
\end{split}
\end{equation}\\
\item\label{Nd}
$(4)$ For $j=0,\cdots,N_{2}^{+}(n,q)$,
\begin{equation*}
\begin{split}
C_{2,j}^{+}(n,q)&=\sum
\binom{q^{-1}A^{+}(n,q)(B^{+}(n,q)+q\delta(2,q;0)+(q-1)^3)}{\nu_{-1},\mu_{-1}}\\
\end{split}
\end{equation*}
\begin{equation}\label{i2}
\begin{split}
&\quad\times \prod_{\beta\neq
-1}\binom{q^{-1}A^{+}(n,q)(B^{+}(n,q)+q\delta(2,q;\beta+1)-2q^2+3q-1)}{\nu_{\beta},\mu_{\beta}}
\end{split}
\end{equation}\\
(cf.(\ref{m1})),\\
where the sum is over all the sets of nonnegative integers
$\{\nu_{\beta}\}_{\beta \in \mathbb {F}_{q}}$ and
$\{\mu_{\beta}\}_{\beta \in \mathbb {F}_{q}}$ satisfying\\
\begin{equation*}
\sum_{\beta\in\mathbb{F}_q} \nu_\beta
+\sum_{\beta\in\mathbb{F}_q}\mu_\beta=j,\quad \text {and} \quad
\sum_{\beta\in\mathbb{F}_q} \nu_\beta
\beta=\sum_{\beta\in\mathbb{F}_q}\mu_\beta \beta.\\
\end{equation*}
\end{theorem}
\begin{theorem}\label{O}$($\cite{D4}$)$
$(1)$ For each odd $n\geq1$ and all $q$, and $h=1,2,3,\cdots,$
\begin{equation}\label{j2}
\begin{split}
2(&\frac{2}{3})^{h}A^{-}(n,q)^{h}\sum^{h}_{j=0}(-1)^{j}\binom{h}{j}B^{-}(n,q)^{h-j}SK^{j}\\
&=q\sum^{min\{N^{-}_{i}(n,q),h\}}_{j=0}(-1)^{j}C^{-}_{j}(n,q)\sum^{h}_{t=j}t!S(h,t)3^{-t}2^{t-j}
\binom{N^{-}_{i}(n,q)-j}{N^{-}_{i}(n,q)-t},
\end{split}
\end{equation}
where  $N_{i}^{-}(n,q)=|DC_{i}^{-}(n,q)|=A^{-}(n,q)B^-(n,q)$,
$S(h,t)$ indicates the Stirling number of the second kind as in
(\ref{n}), and $\{C_{j}^{-}(n,q)\}_{j=0}^{N_{i}^{-}(n,q)}$ denotes
the weight
distribution of $C(DC^{-}(n,q))$ given by\\
\begin{equation*}
\begin{split}
C_{j}^{-}&(n,q)=\sum
\binom{q^{-1}A^{-}(n,q)(B^{-}(n,q)+1)}{\nu_{1},\mu_{1}}\binom{q^{-1}A^{-}(n,q)(B^{-}(n,q)+1)}{\nu_{-1},\mu_{-1}}\\
&\times \prod_{\beta^{2}-1\neq
0~square~}\binom{q^{-1}A^{-}(n,q)(B^{-}(n,q)+q+1)}{\nu_{\beta},\mu_{\beta}}\\
&\times \prod_{\beta^{2}-1
~nonsquare~}\binom{q^{-1}A^{-}(n,q)(B^{-}(n,q)-q+1)}{\nu_{\beta},\mu_{\beta}}~(j=0,\cdots,N_{i}^{-}(n,q)).
\end{split}
\end{equation*}
Here the sum is over all the sets of nonnegative integers
$\{\nu_{\beta}\}_{\beta \in \mathbb {F}_{q}}$ and
$\{\mu_{\beta}\}_{\beta \in \mathbb {F}_{q}}$ satisfying
$\sum_{\beta\in\mathbb{F}_q} \nu_\beta
+\sum_{\beta\in\mathbb{F}_q}\mu_\beta=j$, and
 $\sum_{\beta\in\mathbb{F}_q} \nu_\beta
\beta=\sum_{\beta\in\mathbb{F}_q}\mu_\beta \beta.$\\

\item\label{Ob}
$(2)$ For each even $n\geq2$ and all $q$, and $h=1,2,3,\cdots,$
\begin{equation}\label{k2}
\begin{split}
2&(\frac{2}{3})^{h}A^{+}(n,q)^{h}\sum^{h}_{j=0}(-1)^{j}\binom{h}{j}(B^{+}(n,q)-q^{2}+q)^{h-j}SK^{2j}\\
&=q\sum^{min\{N^{+}_{i}(n,q),h\}}_{j=0}(-1)^{j}C^{+}_{j}(n,q)\sum^{h}_{t=j}t!S(h,t)3^{-t}2^{t-j}
\binom{N^{+}_{i}(n,q)-j}{N^{+}_{i}(n,q)-t},
\end{split}
\end{equation}
where  $N_{i}^{+}(n,q)=|DC_{i}^{+}(n,q)|=A^{+}(n,q)B^+(n,q)$, and
$\{C_{j}^{+}(n,q)\}_{j=0}^{N_{i}^{+}(n,q)}$ is the weight
distribution of $C(DC^{+}(n,q))$ given by
\begin{equation*}
\begin{split}
C_{j}^{+}&(n,q)\\&=\sum
\binom{q^{4}(\delta(2,q;0)+q^5-q^2-3q+3)}{\nu_{0},\mu_{0}}\\
&\times \prod_{\beta \in
\mathbb{F}_{q}^{*}}\binom{q^{4}(\delta(2,q;\beta)+q^5-q^3-q^2-2q+3)}{\nu_{\beta},\mu_{\beta}}~(j=0,\cdots,N_{i}^{+}(n,q)).
\end{split}
\end{equation*}
Here the sum is over all the sets of nonnegative integers
$\{\nu_{\beta}\}_{\beta \in \mathbb {F}_{q}}$ and
$\{\mu_{\beta}\}_{\beta \in \mathbb {F}_{q}}$ satisfying
$\sum_{\beta\in\mathbb{F}_q} \nu_\beta
+\sum_{\beta\in\mathbb{F}_q}\mu_\beta=j$,  and
  $\sum_{\beta\in\mathbb{F}_q} \nu_\beta
\beta=\sum_{\beta\in\mathbb{F}_q}\mu_\beta \beta,$ and, for every
$\beta \in \mathbb{F}_{q}$,
\begin{equation*}
\delta(2,q;\beta)=|\{(\alpha_{1}^{},\alpha_{2}^{})\in
(\mathbb{F}_{q}^{*})^2|\alpha_{1}^{}+\alpha_{1}^{-1}+\alpha_{2}^{}+\alpha_{2}^{-1}=\beta\}|.
\end{equation*}
\end{theorem}

\begin{remark}\label{P}
In \cite{D4},  Theorem 5.4 (1) above is stated to hold for each odd
$n\geq3$ and all $q$, but it is also true for $n=1$ and all $q$.
Indeed, this can be shown by employing the same method as was done
in the proof of Theorem \ref{K} .
\end{remark}

  We are now ready to apply the Pless power moment identity in (\ref{c2}) to\\
$C(DC_{i}^{\mp}(n,q))^{\perp}$, for $i=1,2,$ in order to obtain the
results in Theorem \ref{A}(cf. (\ref{j})-(\ref{m}),
(\ref{o})-(\ref{r})) about  recursive
formulas.\\
  The left hand side of that identity in (\ref{c2}) is equal to
\begin{equation}\label{l2}
\sum_{a\in\mathbb{F}_{q}^{*}}w(c^{\mp}_{i}(a))^{h},
\end{equation}
with the $w(c^{\mp}_{i}(a))$ given by (\ref{d2}) and (\ref{e2}). We
do this for $w(c^{-}_{i}(a))$. In below, $\lq\lq$the sum over
$tra=0$(resp. $tra\neq0$) " will mean $\lq\lq$the sum over all
$a\in\mathbb{F}_{q}^{*}$, with $tra=0$(resp. $tra\neq 0$)."
\begin{equation}\label{m2}
\begin{split}
&\sum_{a\in\mathbb{F}_{q}^{\ast}}w(c^{-}_{i}(a))^{h}=
(\frac{2}{3})^{h}A^{-}(n,q)^{h}\sum_{a\in\mathbb{F}^{\ast}_{q}}(B^{-}(n,q)-(Re\lambda(a))K(\lambda;a^{2}))^{h}\\
&\qquad\quad\quad\quad\qquad=(\frac{2}{3})^{h}A^{-}(n,q)^{h}\sum_{tra=0}(B^{-}(n,q)-K(\lambda;a^{2}))^{h}\\
&\quad\quad\qquad\qquad\qquad+(\frac{2}{3})^{h}A^{-}(n,q)^{h}\sum_{tra\neq0}(B^{-}(n,q)+\frac{1}{2}K(\lambda;a^{2}))^{h}\\
&(\text {nothing that}~ Re\lambda(a)=1, \text ~{if}~ tra=0 ;
Re\lambda(a)=-\frac{1}{2},\text ~{if}~ tra\neq0,i.e.,tra=1,2)\\
&=(\frac{2}{3})^{h}A^{-}(n,q)^{h}\sum_{tra=0}\sum^{h}_{j=0}(-1)^{j}\binom{h}{j}B^{-}(n,q)^{h-j}K(\lambda;a^{2})^{j}\\
&\quad\quad+(\frac{2}{3})^{h}A^{-}(n,q)^{h}\sum_{tra\neq0}\sum^{h}_{j=0}\binom{h}{j}B^{-}(n,q)^{h-j}2^{-j}K(\lambda;a^{2})^{j}\\
&=(\frac{2}{3})^{h}A^{-}(n,q)^{h}\sum_{j=0}^{h}(-1)^{j}\binom{h}{j}B^{-}(n,q)^{h-j}(2SK^{j}-T_{12}SK^{j})
 (cf.(\ref{a}),(\ref{b}))\\
&\quad\quad+(\frac{2}{3})^{h}A^{-}(n,q)^{h}\sum_{j=0}^{h}\binom{h}{j}B^{-}(n,q)^{h-j}2^{-j}T_{12}SK^{j}\\
&=2(\frac{2}{3})^{h}A^{-}(n,q)^{h}\sum_{j=0}^{h}(-1)^{j}\binom{h}{j}B^{-}(n,q)^{h-j}SK^{j}\\
&\quad\quad+(\frac{2}{3})^{h}A^{-}(n,q)^{h}\sum_{j=1}^{h}((-1)^{j+1}+2^{-j})\binom{h}{j}B^{-}(n,q)^{h-j}T_{12}SK^{j}\\
&=q\sum^{min\{N^{-}_{i}(n,q),h\}}_{j=0}(-1)^{j}C^{-}_{j}(n,q)\sum^{h}_{t=j}t!S(h,t)3^{-t}2^{t-j}\binom{N^{-}_{i}(n,q)-j}{N^{-}_{i}(n,q)-t}~(\text{from}~(\ref{j2}))\\
&\quad\quad+(\frac{2}{3})^{h}A^{-}(n,q)^{h}\sum_{j=1}^{h}((-1)^{j+1}+2^{-j})\binom{h}{j}B^{-}(n,q)^{h-j}T_{12}SK^{j}.
\end{split}
\end{equation}
Similarly,
\begin{equation}\label{n2}
\begin{split}
\sum_{a\in\mathbb{F}_{q}^{*}}w(c^{+}_{i}(a)&)^{h}=2(\frac{2}{3})^{h}A^{+}(n,q)^{h}\sum_{j=0}^{h}(-1)^{j}\binom{h}{j}(B^{+}(n,q)-q^{2}+q)^{h-j}SK^{2j}\\
&+(\frac{2}{3})^{h}A^{+}(n,q)^{h}\sum_{j=0}^{h}\binom{h}{j}\{(-1)^{j+1}(B^{+}(n,q)-q^{2}+q)^{h-j}\\
&+2^{-j}(B^{+}(n,q)+\frac{1}{2}q^{2}-\frac{1}{2}q)^{h-j}\}T_{12}SK^{2j}\\
&=q\sum^{min\{N^{+}_{i}(n,q),h\}}_{j=0}(-1)^{j}C^{+}_{j}(n,q)\sum^{h}_{t=j}t!S(h,t)3^{-t}2^{t-j}\binom{N^{+}_{i}(n,q)-j}{N^{+}_{i}(n,q)-t}\\
&+(\frac{2}{3})^{h}A^{+}(n,q)^{h}\sum_{j=0}^{h}\binom{h}{j}\{(-1)^{j+1}(B^{+}(n,q)-q^{2}+q)^{h-j}\\
&+2^{-j}(B^{+}(n,q)+\frac{1}{2}q^{2}-\frac{1}{2}q)^{h-j}\}T_{12}SK^{2j}~(\text{from}~(\ref{k2})).
\end{split}
\end{equation}

On the other hand, the right hand side of  (\ref{c2}) is

\begin{equation}\label{o2}
q\sum^{min\{N^{\mp}_{i}(n,q),h\}}_{j=0}(-1)^{j}C^{\mp}_{i,j}(n,q)\sum^{h}_{t=j}t!S(h,t)3^{-t}2^{t-j}
\binom{N^{\mp}_{i}(n,q)-j}{N^{\mp}_{i}(n,q)-t}.
\end{equation}

Here one has to note that
$dim_{\mathbb{F}_{3}^{}}C(DC_{i}^{\mp}(n,q))=r$ (cf. Theorem
\ref{K}) and to separate the terms corresponding to $j=h$ of the
second sums in (\ref{m2}) and (\ref{n2}). Our main results in
Theorem \ref{A} now follow by equating either (\ref{m2}) or
(\ref{n2}) with (\ref{o2}).

\bibliographystyle{amsplain}

\end{document}